\documentclass[10pt]{amsart}

\usepackage{microtype}
\usepackage[dvipsnames]{xcolor}

\usepackage{amsmath, amssymb, mathrsfs}
\usepackage{mathtools}

\usepackage[mathscr]{euscript}

\newlength{\mylength}
\usepackage{enumitem}
\setlist{listparindent=\parindent, itemsep=0cm, parsep=\mylength, topsep=0cm}

\usepackage[breaklinks=true]{hyperref}
\usepackage{url}

\usepackage{tikz-cd}

\usepackage{amsthm}

\makeatletter
\renewenvironment{proof}[1][\proofname]{\par
	\pushQED{\qed}%
	\normalfont \topsep6\p@\@plus6\p@\relax
	\noindent\emph{#1.}
	\ignorespaces
}{%
\popQED\endtrivlist\@endpefalse
}
\makeatother

\newtheoremstyle{mythm}% name of the style to be used
{\mylength}% measure of space to leave above the theorem. E.g.: 3pt
{0pt}% measure of space to leave below the theorem. E.g.: 3pt
{\itshape}% name of font to use in the body of the theorem
{0pt}% measure of space to indent
{}% name of head font
{. }% punctuation between head and body
{0em}% space after theorem head; " " = normal interword space
{\thmnumber{#2. }\bfseries{\thmname{#1}\thmnote{ (#3)}}}

\newtheoremstyle{mydef}% name of the style to be used
{\mylength}% measure of space to leave above the theorem. E.g.: 3pt
{0pt}% measure of space to leave below the theorem. E.g.: 3pt
{}% name of font to use in the body of the theorem
{0pt}% measure of space to indent
{}% name of head font
{. }% punctuation between head and body
{0em}% space after theorem head; " " = normal interword space
{\thmnumber{#2. }\bfseries{\thmname{#1}\thmnote{ (#3)}}}

\newtheoremstyle{myrmk}% name of the style to be used
{\mylength}% measure of space to leave above the theorem. E.g.: 3pt
{0pt}% measure of space to leave below the theorem. E.g.: 3pt
{}% name of font to use in the body of the theorem
{0pt}% measure of space to indent
{\itshape}% name of head font
{.\ }% punctuation between head and body
{ }% space after theorem head; " " = normal interword space
{\thmname{#1}\thmnumber{ #2}\thmnote{ (#3)}}

\theoremstyle{mythm}
\newtheorem{thm}[subsubsection]{Theorem}
\newtheorem{lem}[subsubsection]{Lemma}
\newtheorem{cor}[subsubsection]{Corollary}

\newtheorem{prop}[subsubsection]{Proposition}
\newtheorem*{thm*}{Theorem}
\newtheorem*{lem*}{Lemma}
\newtheorem*{cor*}{Corollary}
\newtheorem*{claim*}{Claim}
\newtheorem*{prop*}{Proposition}

\theoremstyle{mydef}
\newtheorem{defn}[subsubsection]{Definition}
\newtheorem*{defn*}{Definition}

\theoremstyle{myrmk}

\newtheorem*{rmk*}{Remark}
\newtheorem*{ex*}{Example}

\newcommand{\on}{\operatorname}
\newcommand{\rnc}{\renewcommand}
\rnc{\setminus}{\smallsetminus}

\newcommand{\wt}{\widetilde}

\newcommand{\BB}{\mathbb{B}}
\newcommand{\BZ}{\mathbb{Z}}

\newcommand{\BC}{\mathbb{C}}
\newcommand{\ul}{\underline}

\newcommand{\id}{\on{id}}

\newcommand{\mf}{\mathfrak}
\newcommand{\mc}{\mathscr}

\newcommand{\hra}{\hookrightarrow}

\newcommand{\Hom}{\on{Hom}}

\newcommand{\D}{\mc{D}}

\newcommand{\op}{{\on{op}}}

\DeclareMathOperator*\colim{colim}

\newenvironment{cd}{\begin{equation*}\begin{tikzcd}}{\end{tikzcd}\end{equation*}\ignorespacesafterend}
\newcommand{\e}[1]{\begin{align*} #1 \end{align*}}

\usepackage[margin=1.5in]{geometry}

\makeatletter
\def\blfootnote{\gdef\@thefnmark{}\@footnotetext}
\makeatother

\title[Homotopical presentations of braid groups via reduced lifts]{Homotopical presentations of braid groups \\ via reduced lifts}
\author{James Tao} 
\address{Massachusetts Institute of Technology, Cambridge, MA 02139, USA}
\email{jamestao@mit.edu}
\author{Roman Travkin}
\address{Skolkovo Institute of Science and Technology, Moscow, Russia}
\email{roman.travkin2012@gmail.com}
\date{March 24, 2021}

\definecolor{myblue}{rgb}{0,0.1,0.4}
\newenvironment{myproof}{\color{myblue}\begin{proof}}{\end{proof}}

\newcommand{\nec}{\mathsf{Nec}}

\usepackage{bm}

\newcommand{\subsubsectiona}{\subsubsection{}\hspace{-0.7em}}

\newcommand{\mr}{\mathrm}
\newcommand{\word}{\mathsf{Word}}
\newcommand{\ms}{\mathsf}
\newcommand{\xra}{\xrightarrow} 

\newcommand{\la}{\langle}
\newcommand{\ra}{\rangle}
\newcommand{\acts}{\curvearrowright}

\begin{document}
	
	\begin{abstract}
		In~\cite{deligne2}, Deligne showed that the reduced lift presentation of a finite type generalized braid group remains correct if it is (suitably) interpreted as a presentation of a topological monoid. In this expository paper, we point out that Deligne's argument does not require the `finite type' hypothesis, so it gives a different proof of~\cite[Thm.\ 5.1]{dob}. We also review how to use this result to construct an action of the braid group on the finite or affine Hecke $\infty$-category via intertwining functors. 
	\end{abstract}
	
	\maketitle

	\tableofcontents 
	
	\newpage
	
	\section{Introduction} \label{s-intro} 
	
	\subsection{Reduced lift presentations}  \label{ss-intro1} 
	
	Let $I$ be a Coxeter--Dynkin diagram, which determines a Coxeter group $W_I$, a braid monoid $\BB_I^+$, and a braid group $\BB_I$. Any finite type subset $J \subseteq I$ determines a finite standard subgroup $W_J \subseteq W_I$. The subset of \emph{finite type elements} of $W_I$ is given by  
	\[
		W_{I, \mr{fin}} := \bigcup_{\substack{J \subseteq I \\ \text{$J$ finite type}}} W_J. 
	\]
	The following \emph{reduced lift} presentations of $\BB_I^+$ are well-known: 
	\begin{align}\label{introeq1} 
	\BB_I^+ & \simeq \left\la t_w \text{ for } w \in W_{I, \on{fin}} \ \Bigg| \ \parbox{2.9in}{\begin{center} $t_{w_1} t_{w_2} = t_{w_1w_2}$ for all $w_1, w_2 \in W_{I, \on{fin}}$ satisfying \\ $w_1w_2 \in W_{I, \on{fin}}$ and $\ell(w_1) + \ell(w_2) = \ell(w_1w_2)$ \end{center}} \right\ra  \\ 
	\BB_I^+ & \simeq \left\la t_w \text{ for } w \in W_{I} \ \Bigg| \ \parbox{2.2in}{\begin{center} $t_{w_1} t_{w_2} = t_{w_1w_2}$ for all $w_1, w_2 \in W_{I}$ \\ satisfying $\ell(w_1) + \ell(w_2) = \ell(w_1w_2)$ \end{center}} \right \ra  \label{introeq2} 
	\end{align}
	The right hand sides can be interpreted as presentations of topological monoids, i.e.\ monoid objects in the $\infty$-category $\ms{Spaces}$. \emph{A priori}, the resulting topological monoids may not coincide with the discrete monoid $\BB_I^+$; there is agreement on $\pi_0$, but their higher homotopy groups may not vanish. 
	
	In 1997, Deligne proved that, if $I$ is finite type, then the reduced lift presentation is valid in the $\infty$-category of monoids in $\ms{Spaces}$. This is essentially~\cite[Thm.\ 1.7]{deligne2}, which states that the higher homotopy groups vanish. The main theorem of~\cite{deligne2} is the observation that the 1-truncated version of this statement (pertaining to monoids in 1-groupoids) gives an easy way to construct a monoidal functor from $\BB_I^+$ to any monoidal 1-category. If the language of $\infty$-categories had been available at that time, Deligne could just as well have phrased his main theorem as an easy way to construct a monoidal functor from $\BB_I^+$ to any monoidal $\infty$-category. 
	
	In 2006, Dobrinskaya generalized Deligne's result by removing the assumption that $I$ is finite type. The statement for the presentation~(\ref{introeq1}) is~\cite[Thm.\ 5.1]{dob}, and the statement for presentation~(\ref{introeq2}) is easy to deduce from this. Just as before, one obtains an easy way to construct a monoidal functor from $\BB_I^+$ to any monoidal $\infty$-category. 
	
	The first goal of this paper is to show that Deligne's method from~\cite{deligne2} does not essentially use that $I$ is finite type; this yields an alternative proof of~\cite[Thm.\ 5.1]{dob}. There are pros and cons to both methods. For instance, Deligne's method is more useful when studying factorizations (of an element $b \in \BB_I^+$) with constrained prefixes, as demonstrated in~\cite[Lem.\ 3.6.4]{tatr}. (The present paper is intended to be a convenient reference for~\cite{tatr}.) 
	
	\subsection{Braid group action via intertwining functors}  \label{ss-intro2} 
	
	The results of~\ref{ss-intro1} are useful for constructing monoidal functors from $\BB_I^+$ to any monoidal $\infty$-category. In practice, however, one often wants to construct a monoidal functor from the braid group $\BB_I$. This prompts one to ask the following equivalent questions: 
	\begin{itemize}
		\item Let $\mc{C}$ be any monoidal $\infty$-category. If a monoidal functor $\BB_I^+ \to \mc{C}$ sends each object of $\BB_I^+$ to an invertible object of $\mc{C}$, does it canonically extend to a monoidal functor $\BB_I \to \mc{C}$? 
		\item Is the homotopy groupification\footnote{If $M$ is a topological monoid, then its homotopy groupification can be constructed by taking the classifying $\infty$-category of $M$, inverting all morphisms to obtain an $\infty$-groupoid (i.e.\ topological space), and taking its based loop space. The homotopy groupification is the universal group object in the $\infty$-category $\ms{Spaces}$ which receives a monoid map from $M$.} of the discrete monoid $\BB_I^+$ equivalent to $\BB_I$? 
	\end{itemize}
	When $I$ is finite type, this question was essentially answered affirmatively in~\cite[1.10]{deligne2}, using results of the 1972 paper~\cite{deligne}. Hence, the main theorem of~\cite{deligne2} applies not only to $\BB_I^+$ but also to $\BB_I$. 
	
	What about when $I$ is not finite type? The main theorem of~\cite{dob} states that this question is equivalent to the $K(\pi, 1)$ conjecture, which predicts that the $W_I$-quotient of the complement of a complexified hyperplane arrangement is a classifying space for $\BB_I$. Hence, if the $K(\pi, 1)$ conjecture is known for a particular Coxeter--Dynkin diagram $I$, then one obtains an easy way to construct monoidal functors from $\BB_I$. The main result of~\cite{deligne} was the proof of the $K(\pi, 1)$ conjecture when $I$ is finite type. In 2020, the $K(\pi, 1)$ conjecture was proved when $I$ is affine type, see~\cite{ps}. 
	
	Assume that $I$ is finite or affine type. The second goal of this paper is to use the idea of the previous paragraph to construct a monoidal functor from $\BB_I$ to a (finite or affine) Hecke $\infty$-category with cells indexed by $W_I$, see Corollary~\ref{braid-group}. This amounts to rewriting the proof of `Application 2' in~\cite[\S 0]{deligne2} in $\infty$-categorical language, and no new difficulties arise. 
	
	\subsubsection{Remark} \label{history} 
	This monoidal functor plays an important role in geometric representation theory, and we survey some of the previous literature here. Assume that $I$ is finite type. In 1983, Beilinson and Bernstein defined \emph{intertwining functors} on $\D(G/B)$ and studied their effect on twistings and support of $\D$-modules~\cite[\S 7]{intertwine}. Roughly speaking, an intertwining functor is given by pull-push along a span 
	\[
	G/B \leftarrow X_w \rightarrow G/B
	\]
	where $w \in W$ is any element, and $X_w$ parameterizes pairs of Borel subgroups in relative position $w$. From their construction, one can easily obtain an action of $\BB_I$ on $K_0(\D(G/B))$. In 1995, Brou\'e and Michel~\cite{brmi} used the reduced lift presentation of $\BB_I^+$ to generalize $X_w$ by allowing $w$ to be an element of $\BB_I^+$, see~\cite[1.B]{brmi}. The resulting varieties $X_b$ for $b \in \BB_I^+$ were defined up to an isomorphism. (They also gave an analogous generalization of Deligne--Lusztig varieties, which were the focus of their paper.) Deligne's 1997 paper~\cite{deligne2} showed that the varieties $X_b$ are defined up to a \emph{unique} isomorphism, and the resulting intertwining functors give an action of $\BB_I$ on $\D(G/B)$. This insight was incorporated into the 1997 published version of Brou\'e and Michel's paper. 
	
	In 2004, Rouquier considered the full subcategory $\D^b_\sigma(G/B) \subset \D(G/B)$ of complexes of $\D$-modules whose cohomology sheaves are bounded and constant on the $B$-orbits, showed that the action $\BB_I \acts \D^b_\sigma(G/B)$ can be improved to an action by a `categorified braid group' $\mc{B}_{W_I}$ (defined in~\cite[3.3.2]{rouquier}), and conjectured that the same is true for $\D(G/B)$ in place of $\D^b_\sigma(G/B)$. The object $\mc{B}_{W_I}$ is richer than $\BB_I$ because it is a monoidal 1-category which is not a groupoid; it is defined as a full subcategory of a Soergel bimodule 1-category. We will not discuss $\mc{B}_{W_I}$ in the present paper. 
	
	\subsection{Notations and basic definitions} \label{ss-notation} 
	A \emph{Coxeter--Dynkin diagram} is a pair $(I, m)$ where $I$ is a finite set and $m : I \times I \to \BZ \sqcup \{\infty\}$ is a symmetric map such that $m_{ss} = 1$ and $m_{st} \ge 2$ for all distinct $s, t\in I$. For the entire paper, we fix a Coxeter--Dynkin diagram, which we denote as $I$ rather than $(I, m)$.  	
	
	For any Coxeter--Dynkin diagram $I$, there is the associated \emph{Coxeter group} 
	\[
	W_I = \la s \in I \ | \ (st)^{m_{st}} = 1 \text{ for } s, t \in I \ra.
	\]
	There is also the \emph{braid group} (a.k.a.\ Artin group)
	\[
	\BB_I = \la s \in I \ | \ \underbrace{(sts \cdots)}_{m_{st} \text{ factors}} = \underbrace{(tst \cdots)}_{m_{st} \text{ factors}} \text{ for } s, t \in I \ra
	\]
	and the \emph{braid monoid} $\BB_I^+$ which is specified by the same presentation.
	
	The \emph{length} of an element $w \in W_I$, denoted $\ell(w)$, is the smallest possible length of a product of `simple reflection' generators which equals $w$: 
	\begin{equation*} 
	w = s_1 \cdots s_n
	\end{equation*}
	The sequence $(s_1, \ldots, s_n)$ is called a \emph{minimal-length expression} for $w$. Similarly for $\BB_I^+$. 		
	
	Any $J \subseteq I$ defines a subdiagram $(J, m_J)$ with Coxeter group $W_J \subseteq W_I$, braid group $\BB_J \subseteq \BB_I$, and braid monoid $\BB_J^+ \subseteq \BB_I^+$. If $W_J$ is finite, we say that $J$ is \emph{finite type}. We also define subsets of $W_I$ and $\BB_I^+$ given by 
	\e{
		W_{I, \mr{fin}} &:= \bigcup_{\substack{J \subseteq I \\ J \text{ finite type}}}W_J \\
		\BB_{I, \mr{fin}}^{+} &:= \bigcup_{\substack{J \subseteq I \\ J \text{ finite type}}} \BB_J^+
	}
	and call their elements \emph{finite type}. We define an element of $\BB_I^+$ to be \emph{reduced} if it cannot be written as $b_1s_is_ib_2$ for $b_1, b_2 \in \BB_I^+$ and $i \in I$.\footnote{In~\cite[Sect.\ 4]{dob}, such an element is called \emph{squarefree}.} 
	
	The \emph{prefix order} (a.k.a.\ weak right Bruhat order) on $W_I$ is defined by $w' \preceq w$ if there is a minimal-length expression for $w$ (as above) such that some initial substring $s_1 \cdots s_j$ is an expression for $w'$. It is equivalent to require that  
	\[
	\ell(w) = \ell(w') + \ell((w')^{-1}w). 
	\]
	We say that $w'$ is a \emph{prefix} of $w$. There is also a prefix order on $\BB_I^+$. We emphasize that the symbols $\preceq$ and $\prec$ will always refer to the prefix order, not the (strong) Bruhat order. 
	
	There is a well-defined map of sets
	\[
	r : W_I \to \BB_I^+
	\]
	given as follows: given any minimal-length expression for $w \in W_I$, we define $r(w)$ to be the element of $\BB_I^+$ specified by the same expression. This is called the \emph{reduced lift} map because it is a bijection onto the set of reduced elements of $\BB_I^+$. A sequence $(w_1, \ldots, w_n)$ in $W_I$ is called \emph{reduced} if $r(w_1) \cdots r(w_n) \in \BB_I^+$ is reduced.\footnote{This definition is usually only applied to sequences of simple reflections, but the more general scope will be convenient for stating Definition~\ref{dob-present}.} 
	
	\subsection{Acknowledgments} We would like to thank Roman Bezrukavnikov for helpful conversations and for pointing out many of the references surveyed in~\ref{history}. The first author is supported by the NSF GRFP, grant no.\ 1122374. 
	
	\section{Contractibility of factorization posets} \label{s-contract} 
	
	In this section, we establish some analogues of~\cite[Thm.\ 1.7]{deligne2} using the same method of proof. The proof in~\cite[\S 2]{deligne2} is phrased in the language of chambers and galleries, but we will not use this language. 
	
	\subsection{Maximal finite type prefixes} \label{ss-prefix} 
	
	\begin{defn} \label{prefix-def} 
		For $b \in \BB_I^+$, we define the sets 
		\e{
			\mc{L}(b) &:= \{s \in I \, | \, s \preceq  b\} \\
			\ms{Pre}(b) &:= \{b' \in \BB_I^+ \, | \, b' \preceq b\}
		} 
	\end{defn}
	
	\begin{lem} \label{prefix1} 
		Let $b \in \BB_I^+$ be reduced. Then $\mc{L}(b)$ is finite type, and $r(W_{\mc{L}(b)}) \subseteq \ms{Pre}(b)$. 
	\end{lem}
	\begin{myproof}
		This is a restatement of Lemma~4.7.2 and Lemma~4.7.3 in~\cite{davis-book}. 
	\end{myproof}
	
	\begin{cor} \label{prefix2} 
		Let $b \in \BB_I^+$ be arbitrary. Then $\mc{L}(b)$ is finite type, and $r(W_{\mc{L}(b)}) \subseteq \ms{Pre}(b)$. 
	\end{cor}
	\begin{myproof}
		By~\cite[Prop.\ 2.1]{m}, $b$ has a unique maximal reduced prefix, say $b' \in \BB_I^+$. Thus $\mc{L}(b') = \mc{L}(b)$ and $\ms{Pre}(b') \subseteq \ms{Pre}(b)$, so the claim follows from Lemma~\ref{prefix1} applied to $b'$. 
	\end{myproof}
	
	\subsection{Contractibility results} \label{ss-results} 
	We will show that several versions of the `poset of factorizations of an element $b \in \BB_I^+$' are contractible. For conciseness, instead of writing four nearly identical proofs, we use the trick from the proof of~\cite[Prop.\ 5.8]{dob} to deduce several contractibility results from a single result (Theorem~\ref{results2}). 
	
	\begin{defn} \label{results-def} 
		For $b \in \BB_I^+$, let $\ms{Word}(b)$ be the poset of strictly increasing chains 
		\[
			(1 \prec b_1 \prec \cdots \prec b_{n-1} \prec b) 
		\]
		in $\BB_I^+$ which start at $1$ and end at $b$. The partial order on $\ms{Word}(b)$ is the opposite of that given by containment. We also introduce the following full subposets: 
		\begin{itemize}
			\item $\ms{Word}_{\mr{s}}(b)$ is obtained from $\ms{Word}(b)$ by deleting the terminal object $(1 \prec b)$. 
			\item $\ms{Word}_{\mr{r}}(b)$ consists of chains such that each $b_{i-1}^{-1}b_i$ is reduced. 
			\item $\ms{Word}_{\mr{f}}(b)$ consists of chains such that each $b_{i-1}^{-1}b_i$ is finite type. 
			\item $\ms{Word}_{\mr{fr}}(b) := \ms{Word}_{\mr{r}}(b) \cap \ms{Word}_{\mr{f}}(b)$. 
			\item $\ms{Word}_{\Delta}(b)$ consists of chains such that each $b_{i-1}^{-1}b_i$ is $r(\Delta_J)$ for some finite type $J \subseteq I$. 
		\end{itemize}
	\end{defn}
	
	\subsubsection{Remark} \label{results-rmk} 
		As a 1-category, $\ms{Word}(b)$ can be equivalently defined as follows: 
		\begin{itemize}
			\item Objects are sequences $\bm{b} = (b_1, \ldots, b_n)$ in $\BB_I^+ \setminus \{1\}$ which multiply to $b$. 
			\item A morphism $\varphi : (b_{1,1}, \ldots, b_{1, n_1}) \to (b_{2, 1}, \ldots, b_{2, n_2})$ is a weakly increasing surjective map $\varphi_* : [n_1] \to [n_2]$ such that, for each $j \in [n_2]$, we have 
			\[
			b_{2, j} = (\text{product of $b_{1, i}$ for $i \in \varphi_*^{-1}(j)$}). 
			\] 
		\end{itemize}
		The equivalence between these two definitions sends the strictly increasing chain in Definition~\ref{results-def}
		to the sequence $(b_1, b_1^{-1} b_2, \ldots, b_{n-1}^{-1} b)$. 
		
	\begin{prop} \label{results1} 
		If $b \in \BB_I^+ \setminus \{r(\Delta_J) \, | \, \textnormal{finite type } J \subseteq I \}$, then  $\ms{Word}_{\mr{s}}(b)$ is contractible. 
	\end{prop}
	Compare the following proof with that of~\cite[Thm.\ 2.4]{deligne2}. 
	\begin{myproof}
		Consider the functor 
		\[
			F : \ms{Word}_{\mr{s}}(b) \to (\text{nonempty subsets of $\mc{L}(b)$})
		\]
		which sends a chain 
		\[
			(1 \prec b_1 \prec \cdots \prec b_{n-1} \prec b) 
		\]
		to the nonempty subset $\mc{L}(b_1) \subseteq \mc{L}(b)$. We will show that, for each nonempty subset $T \subseteq \mc{L}(b)$, the slice poset $(T \downarrow F)$ is contractible. Then Quillen's Theorem A implies that $F$ is a homotopy equivalence. Since the target of $F$ has a terminal object, it is contractible. This proves that $\ms{Word}_{\mr{s}}(b)$ is contractible. 
		
		Now we begin the proof that $(T \downarrow F)$ is contractible. Let $\mc{C} \subseteq (T \downarrow F)$ be the full subposet consisting of chains for which $b_1 = r(\Delta_T)$. There is an adjunction 
		\[
			\mc{C} \rightleftarrows (T \downarrow F)
		\]
		where the left adjoint is the embedding and the right adjoint modifies a chain via 
		\[
			(1 \prec b_1 \prec \cdots ) \mapsto \begin{cases}
				(1 \prec r(\Delta_T) \prec b_1 \prec \cdots) & \text{if } b_1 \neq r(\Delta_T) \\
				(1 \prec b_1 \prec \cdots ) & \text{if } b_1 = r(\Delta_T)
			\end{cases} 
		\]
		The right adjoint is well-defined by Corollary~\ref{prefix2}, which shows that $r(\Delta_T) \prec b_1$. Quillen's Theorem A implies that these adjunctions are homotopy equivalences. Finally, note that $\mc{C}$ has a terminal object given by 
		\[
			(1 \prec r(\Delta_T) \prec b). 
		\]
		(The hypothesis that $b \neq r(\Delta_J)$ for any finite type $J \subseteq I$ is used to show that this chain lies in $\ms{Word}_{\mr{s}}(b)$.) Therefore, $\mc{C}$ is contractible, and so is $(T \downarrow F)$. 		
	\end{myproof}
	
	\begin{thm} \label{results2} 
		Let $b \in \BB_I^+$ be arbitrary. Any full subposet of $\ms{Word}(b)$ which contains $\ms{Word}_{\Delta}(b)$ is contractible. 
	\end{thm}
	\begin{myproof}
		Let $\bm{x}_1, \ldots, \bm{x}_n \in \ms{Word}(b) \setminus \ms{Word}_{\Delta}(b)$ be the objects which are not contained in the given subposet, sorted so that their sizes weakly increase. For each $i$, let $\word_{}(b)_{\ge i} \subset \word_{}(b)$ be the full subposet obtained by deleting $\bm{x}_1, \ldots, \bm{x}_{i-1}$. We will show that each embedding
		\[
		F : \word_{}(b)_{{\ge}i+1} \hra \word_{}(b)_{\ge i}
		\]		
		is a homotopy equivalence.
		
		By Quillen's Theorem A, it suffices to show that $(F \downarrow \bm{x}_i)$ is contractible. Since the $\bm{x}_i$ are listed in order, we have 
		\[
			(F \downarrow \bm{x}_i) \simeq (\word_{}(b) \downarrow \bm{x}_i) \setminus \{\bm{x}_i\}.
		\]
		If $\bm{x}_i = (1 = x_{i, 0} \prec x_{i, 1} \prec \cdots \prec x_{i, k-1} \prec x_{i, k} = b)$, then $(F \downarrow \bm{x}_i)$ is the full subposet of
		\[
			\word(x_{i, 1}) \times \word(x_{i, 1}^{-1} x_{i, 2}) \times \cdots \times \word(x_{i, k-1}^{-1} b)
		\]
		consisting of tuples such that at least one coordinate is a chain of size $\ge 3$. Let $S \subseteq \{1, \ldots, k\}$ be the index set such that $j \in S$ if and only if $x_{i, j-1}^{-1}x_{i, j} \notin \{r(\Delta_J) \, | \, \textnormal{finite type } J \subseteq I \}$. Since $\bm{x}_i \notin \word_{\Delta}(b)$, the set $S$ is nonempty. For $j = 1, \ldots, k$, define the poset $\mc{C}_j$ as follows:
		\begin{enumerate}[label=(\roman*)]
			\item If $j\notin S$, then $\mc{C}_j = \word(x_{i, j-1}^{-1}x_{i, j})$.
			\item If $j \in S$, then $\mc{C}_j = \word_{\mr{s}}(x_{i, j-1}^{-1}x_{i, j})$.
		\end{enumerate}
		Note that each $\mc{C}_j$ is contractible. For (i) this is because $\word(x_{i, j-1}^{-1}x_{i, j})$ has a terminal object, and for (ii) this is Proposition~\ref{results1}. Consider the full embedding
		\[
		G: \underset{j}{\times} \, \mc{C}_j \hra (F \downarrow \bm{x}_i).
		\]
		We will show that $G$ is a homotopy equivalence.
		
		Let $\mf{y} = (\bm{y}_1, \ldots, \bm{y}_k) \in (F \downarrow \bm{x}_i)$ be any object not contained in the image of $G$, where the notation means that $\bm{y}_j \in \ms{Word}(x_{i, j-1}^{-1} x_{i,j})$ for each $j$. We have 
		\[
			(G \downarrow \mf{y}) \simeq \underset{j}{\times}\, \mc{D}_j
		\]
		where $\mc{D}_j$ for $j = 1, \ldots, k$ is described as follows:	
		\begin{enumerate}[label=(\roman*)] \setcounter{enumi}{2}
			\item If $\bm{y}_j$ has size 2, then $\mc{D}_j = \mc{C}_j$.
			\item If $\bm{y}_j$ has size $\ge 3$, then $\mc{D}_j$ has a terminal object.
		\end{enumerate}
		In either case, $\mc{D}_j$ is contractible. Therefore $(G \downarrow \mf{y})$ is contractible, and Quillen's Theorem A implies that $G$ is a homotopy equivalence.
		
		Since the $\mc{C}_j$ are contractible, so is $(F \downarrow \bm{x}_i)$. Therefore, Quillen's Theorem A implies that $F$ is a homotopy equivalence, as desired. 		
	\end{myproof}
	
	\begin{cor} \label{results3} 
		Let $b \in \BB_I^+$ be arbitrary. Then $\ms{Word}_{\mr{r}}(b), \ms{Word}_{\mr{f}}(b),  \ms{Word}_{\mr{fr}}(b),  \ms{Word}_{\Delta}(b)$ are contractible. 
	\end{cor}
	\begin{myproof}
		These are full subposets of $\ms{Word}(b)$ which contain $\ms{Word}_{\Delta}(b)$, so the claim follows from Theorem~\ref{results2}. 
	\end{myproof}
	
	\section{The braid group as a homotopy colimit} \label{s-colimit} 
	
	In~\ref{ss-neck} and \ref{ss-partial}, we give a fairly self-contained account of the homotopy monoidification of a partial monoid; for the history of this construction, which goes back to Segal, see the discussion in~\cite{dob}. Using this framework, we show that the contractibility results of~\ref{ss-results} imply that the reduced lift presentation of $\BB_I^+$ is correct when (suitably) interpreted as a presentation of a homotopy monoid. This follows the proof of~\cite[Thm.\ 1.5]{deligne2}. 
	
	We also record two applications of this result: a presentation of $\BB_I^+$ as the homotopy colimit of its finite type standard submonoids (Corollary~\ref{cors1}) and a monoidal functor from $\BB_I$ to the Hecke $\infty$-category via intertwining functors when $I$ is finite or affine type (Corollary~\ref{braid-group}). 
	
	\subsection{Fibrant replacement via necklaces} \label{ss-neck} 
	This subsection reviews a method for computing mapping spaces in the Joyal fibrant replacement of a simplicial set. 
	
	\begin{defn} \label{neck-def} 
		A \emph{necklace} is a nonempty simplicial set of the form
		\[
		\Delta^{n_0} \vee \Delta^{n_1} \vee \cdots \vee \Delta^{n_k}
		\]
		where each wedge means that the vertex $n_i \in \Delta^{n_i}$ is identified with the vertex $0 \in \Delta^{n_{i+1}}$. Each necklace has a unique `first vertex' $0 \in \Delta^{n_0}$ and `last vertex' $n_k \in \Delta^{n_k}$. The simplicial set $\Delta^0$ is the unique necklace whose first and last vertices coincide. 
		
		Let $\nec$ be the 1-category whose objects are necklaces and whose morphisms are maps of simplicial sets which send the first vertex to the first vertex and the last vertex to the last vertex. If $S$ is a simplicial set, and $a, b \in S_0$ are vertices, we define $\nec_S(a, b)$ to be the 1-category of pairs $(T, \pi)$ where $T \in \nec$ and $\pi : T \to S$ is a map sending the first vertex to $a$ and the last vertex to $b$. As usual, morphisms in $\nec_S(a, b)$ are commutative triangles. 
	\end{defn}
	
	\begin{thm} \label{neck1} 
		Let $S$ be a simplicial set, let $a, b \in S_0$ be vertices, and let $S \to \wt{S}$ be a fibrant replacement in the Joyal model structure. Then $\Hom_{\wt{S}}(a, b)$ is homotopy equivalent to $\nec_S(a, b)$. 
	\end{thm}
	\begin{myproof}
		Theorem~1.3 in~\cite{rigid} states that there is a zig-zag of weak equivalences between two simplicial categories $\mf{C}^{\mr{nec}}(S)$ and $\mf{C}[S]$. The simplicial category $\mf{C}^{\mr{nec}}(S)$ is defined so that 
		\[
		\Hom_{\mf{C}^{\mr{nec}}(S)}(a, b) = \on{N}(\nec_S(a, b)),
		\]
		while the simplicial category $\mf{C}[S]$ is defined in~\cite[1.1.5]{htt} and has the property that the topological nerve of $\lvert\mf{C}[S]\rvert$ is a fibrant replacement for $S$. 
	\end{myproof}
	
	\begin{defn} \label{neck2} 
		For any simplex $\Delta^n$, there is the necklace 
		\[
		\on{spine}(\Delta^n) := \Delta^{\{0, 1\}} \vee \Delta^{\{1, 2\}} \vee \cdots \vee \Delta^{\{n-1, n\}} \subset \Delta^n
		\]
		which we call a \emph{spine}. We introduce a property which applies to a simplicial set $S$:  
		\begin{enumerate}
			\item[(P0)] For any $n \ge 2$, in any solid diagram as shown below, there is at most one dashed map which makes the diagram commute:\footnote{This condition is vacuously true if $n = 0, 1$.}  
			\begin{cd}
				\on{spine}(\Delta^n) \ar[d, hookrightarrow] \ar[r] & S \\
				\Delta^n \ar[ru, dashed] 
			\end{cd}
		\end{enumerate}
		If $\varphi : \on{spine}(\Delta^n) \to S$ is a map for which the above lift exists, we call $\varphi$ a \emph{strongly composable} sequence of edges of $S$, and $\Delta^{\{0, n\}} \to S$ is their \emph{composite}. 
	\end{defn}
	
	\subsubsectiona \label{neck3} 
	Let $S$ be a simplicial set which satisfies (P0), and let $a, b\in S_0$ be vertices. Then $\ms{Nec}_S(a, b)$ is equivalent to the following 1-category: 
	\begin{itemize} 
		\item An object is a triple $(\Delta^n, \varphi, p)$ where 
		\[
		\on{spine}(\Delta^n) \xra{\varphi} S
		\]
		is a map satisfying $\varphi(0) = a$ and $\varphi(n) = b$, and $p$ is an ordered partition
		\[
		n = n_{p, 1} + n_{p, 2} + \cdots + n_{p, k}
		\]
		into positive integers $n_{p, i} \ge 1$, satisfying that, for each $i = 1, \ldots, k$, the subsequence of arrows 
		\[
		\on{spine}(\Delta^{n_{p, i}}) \xra{\varphi} S
		\]
		joining the vertices $n_{p, 1} + \cdots + n_{p, i-1}$ and $n_{p,1} + \cdots + n_{p, i}$ in $\Delta^n$ 	is strongly composable. 
		\item A morphism $F : (\Delta^{n_1}, \varphi_1, p_1) \to (\Delta^{n_2}, \varphi_2, p_2)$ is a map $F_* : \Delta^{n_1} \to \Delta^{n_2}$ satisfying the following properties: 
		\begin{enumerate}[label=(\roman*)]
			\item The partition $p_1$ pushes forward under $F_*$ in a natural way, yielding a partition $F_*(p_1)$ of $n_2$. We require that $F_*(p_1)$ refines $p_2$. 
			\item For any integer $0 \le j < n_1$, (i) implies that the sequence of arrows 
			\[
			\on{spine}(\Delta^{\{F_*(j), F_*(j)+1, \ldots, F_*(j+1)\}}) \xra{\varphi_2} S
			\]
			is strongly composable. We require that its composite equals $\varphi_1(\Delta^{\{j, j+1\}})$.  
		\end{enumerate}
	\end{itemize} 
	To see that these 1-categories are equivalent, note that, if 
	\[
	\on{spine}(\Delta^n) \xra{\varphi} S 
	\]
	is a map, and $p$ is an ordered partition of $n$ such that the arrows $\varphi(\on{spine}(\Delta^{n_{p, i}}))$ are composable, then $\varphi$ extends uniquely to a map $\wt{\varphi}$ as follows: 
	\[
	\Delta^{n_{p, 1}} \vee \Delta^{n_{p, 2}} \vee \cdots \vee \Delta^{n_{p, k}} \xra{\wt{\varphi}} S. 
	\]

	\subsection{Fibrant replacement of partial 1-categories} \label{ss-partial} 
	We consider a special case of the situation described in~\ref{ss-neck}. If a simplicial set is a \emph{partial 1-category}, then the mapping spaces in its Joyal fibrant replacement are especially easy to compute. If the simplicial has only one vertex, this recovers the notion of `homotopy monoidification' of a partial monoid. 
	
	\begin{defn} \label{partial-def} 
		A simplicial set $C$ is a \emph{partial 1-category} if it satisfies the following: 
		\begin{enumerate}[label=(P\arabic*)] 
			\item In any solid diagram as shown below, there is at most one dashed map which makes the diagram commute: 
			\begin{cd}
				\Lambda_1^2 \ar[d, hookrightarrow] \ar[r] & C \\
				\Delta^2 \ar[ru, dashed] 
			\end{cd}
			\item In any solid diagram as shown below, there is exactly one dashed map which makes the diagram commute: 
			\[
			\begin{tikzcd}
			\Delta^{\{0, 1, 2\}} \underset{\Delta^{\{0, 2\}}}{\cup} \Delta^{\{0, 2, 3\}} \ar[d, hookrightarrow] \ar[r] & C \\
			\Delta^3 \ar[ru, dashed] 
			\end{tikzcd}
			\qquad 
			\begin{tikzcd}
			\Delta^{\{0, 1, 3\}} \underset{\Delta^{\{1, 3\}}}{\cup} \Delta^{\{1, 2, 3\}} \ar[d, hookrightarrow] \ar[r] & C \\
			\Delta^3 \ar[ru, dashed] 
			\end{tikzcd}
			\]
			\item $C$ is 2-cotruncated, i.e.\ for any $n \ge 3$, in the following solid diagram, there is exactly one dashed map making the diagram commute: 
			\begin{cd}
				\partial \Delta^n \ar[d, hookrightarrow] \ar[r] & C \\
				\Delta^n \ar[ru, dashed] 
			\end{cd}
		\end{enumerate}
		Note that (P1) is a special case of (P0). In (P2), the upper-left objects look like 
		\[ 
		\begin{tikzcd}
		1 \ar[r] & 2 \ar[d] \\
		0 \ar[u] \ar[ru] \ar[r]  & 3
		\end{tikzcd}
		\qquad \text{ and } \qquad
		\begin{tikzcd}
		1 \ar[r] \ar[rd] & 2 \ar[d] \\
		0 \ar[u] \ar[r] & 3
		\end{tikzcd}
		\]
		These should be thought of as two kinds of `commutative squares' in a simplicial set. 
	\end{defn}
	
	\subsubsectiona \label{partial1} 
	Let $C$ be a partial 1-category, and consider a sequence of arrows 
	\[
		\varphi : \mr{spine}(\Delta^n) \to C. 
	\]
	A full parenthesization of the expression 
	\[
		\varphi((n-1) \to n) \circ \varphi((n-2) \to (n-1)) \circ \cdots \circ \varphi(0 \to 1)
	\]
	corresponds to a triangulation of an $(n+1)$-gon $P_{n+1}$ with vertices numbered $0, 1, \ldots, n$. We say that the \emph{parenthesized composite} of these arrows exists in $C$ if there is a lift 
	\begin{cd}
		\mr{spine}(\Delta^n) \ar[r, "\varphi"] \ar[d, hookrightarrow] & C \\
		P_{n+1} \ar[ru, dashed] 
	\end{cd}
	where $P_{n+1}$ is viewed as a simplicial set via the aforementioned triangulation. If this lift exists, it is unique by (P1). 
	
	Moreover, (P2) implies that, if the parenthesized composite exists for one parenthesization, then it exists for all parenthesizations, and the resulting arrows in $C$ are equal. This uses the standard argument for bootstrapping from a 3-term associativity law to an $n$-term associativity law. When this is the case, (P3) implies that the arrows in $\varphi$ are strongly composable (Definition~\ref{neck2}), and that the lift 
	\begin{cd}
		\mr{spine}(\Delta^n) \ar[r] \ar[d, hookrightarrow] & C \\
		\Delta^n \ar[ru, dashed] 
	\end{cd}
	is unique. In particular, $C$ satisfies (P0). Conversely, if the arrows in $\varphi$ are strongly composable, then their parenthesized composite exists. 
	
	\begin{defn} \label{spine} 
		If $C$ is a partial 1-category, and $a, b \in C_0$ are vertices, we define $\ms{Spine}_C(a, b)$ to be a 1-category given as follows: 
		\begin{itemize}
			\item The objects are pairs $(\Delta^n, \varphi)$ where $\varphi : \on{spine}(\Delta^n) \to C$ is a map sending the first vertex to $a$ and the last vertex to $b$. 
			\item A morphism $F : (\Delta^{n_1}, \varphi_1) \to (\Delta^{n_2}, \varphi_2)$ is a map $F^* : \Delta^{n_2} \to \Delta^{n_1}$ satisfying the following property: 
			\begin{itemize}
				\item For any integer $0 \le j < n_1$, the composite of the arrows 
				\[
				\varphi_1(\Delta^{\{F^*(j), F^*(j)+1, \ldots, F^*(j+1)\}}) \to C
				\]
				exists and equals $\varphi_2(\Delta^{\{j, j+1\}})$.
			\end{itemize}
			\item For two morphisms 
			\begin{cd}
				(\Delta^{n_1}, \varphi_1) \ar[r, "F_{21}"] & (\Delta^{n_2}, \varphi_2) \ar[r, "F_{32}"] & (\Delta^{n_3}, \varphi_3)
			\end{cd}
			the composite $F_{32} \circ F_{21}$ is given by the map $F_{21}^* \circ F_{32}^* : \Delta^{n_3} \to \Delta^{n_1}$. This map satisfies the previous condition labeled by `--' because the composite of a sequence of morphisms in $C$ can be computed using any parenthesization, see~\ref{partial1}. 
		\end{itemize} 
		There is a functor $U : \ms{Nec}_C(a, b) \to \ms{Spine}_C(a, b)^{\op}$ defined as follows: when $\ms{Nec}_S(a, b)$ is defined as in~\ref{neck3}, the functor $U$ sends $(\Delta^n, \varphi, p) \mapsto (\Delta^n, \varphi)$. 
	\end{defn} 
	
	\begin{lem} \label{partial2} 
		Let $C$ be a partial 1-category, and fix vertices $a, b \in C_0$. The functor $U : \ms{Nec}_C(a, b) \to \ms{Spine}_C(a, b)^{\op}$ is a homotopy equivalence. 
	\end{lem}
	\begin{myproof}
		Redefine $\ms{Nec}_S(a, b)$ as in~\ref{neck3}. First, we show that $U$ is cocartesian. Given an object $(\Delta^{n_1}, \varphi_1, p_1) \in \ms{Nec}_C(a, b)$ and a morphism $(\Delta^{n_1}, \varphi_1) \to (\Delta^{n_2}, \varphi_2)$ in $\ms{Spine}_C(a, b)^{\op}$, which is encoded by a map $F_* : \Delta^{n_1} \to \Delta^{n_2}$, a cocartesian lift of this morphism is given by
		\[
		(\Delta^{n_1}, \varphi_1, p_1) \to (\Delta^{n_2}, \varphi_2, F_*(p_1))
		\]
		where $F_*(p_1)$ is the pushforward of $p_1$, as in~\ref{neck3}. 
		
		Since $U$ is cocartesian, Quillen's Theorem A implies that, if each fiber of $U$ is contractible, then $U$ is a homotopy equivalence. To finish, we will show that each fiber of $U$ is contractible. For any $(\Delta^n, \varphi) \in \ms{Spine}_C(a, b)$, the fiber category $U^{-1}((\Delta^n, \varphi))$ is the poset of partitions $p$ for which each sequence of arrows $\on{spine}(\Delta^{n_{p_i}})$ is composable in $C$, equipped with the usual partial order on partitions. This poset is contractible because it has an initial object, given by the partition of $n$ into parts of size 1. 		
	\end{myproof}
	
	\begin{cor} \label{partial3} 
		Let $C$ be a partial 1-category, let $a, b \in C_0$ be vertices, and let $C \to \wt{C}$ be a fibrant replacement in the Joyal model structure. Then there is a homotopy equivalence 
		\[
		\Hom_{\wt{C}}(a, b) \simeq \on{N}(\ms{Spine}_C(a, b)). 
		\]
	\end{cor}
	\begin{myproof}
		Follows from Theorem~\ref{neck1} and Lemma~\ref{partial3}. 
	\end{myproof}
	
	\subsubsectiona \label{partial4} 
	Let $C$ be a partial 1-category, and fix vertices $a, b \in C_0$. Define the full subcategory 
	\[
	\ms{Spine}_C^{\mr{nd}}(a, b) \subset \ms{Spine}_C(a, b)
	\]
	to consists of maps $\varphi : \on{spine}(\Delta^n) \to C$ which send each edge to a nondegenerate edge of $C$. 
	\begin{lem*} 
		Suppose there is no composable pair of nondegenerate edges of $C$ whose composite is degenerate. If $a \neq b$, we have a homotopy equivalence 
		\[
		\ms{Spine}_C^{\mr{nd}}(a, b) \to \ms{Spine}_C(a, b). 
		\]
		If $a = b$, we have a homotopy equivalence 
		\[
		\ms{Spine}_C^{\mr{nd}}(a, a)  \sqcup \{*\} \to \ms{Spine}_C(a, a)
		\]
		which sends $*$ to the unique degenerate edge at $a$. 		
	\end{lem*}
	\begin{myproof}
		Assume that $a \neq b$. There is an adjunction
		\[
		\ms{Spine}_C^{\mr{nd}}(a, b) \rightleftarrows \ms{Spine}_C(a, b). 
		\]
		where the left adjoint is the embedding, and the right adjoint sends a map $\varphi : \on{spine}(\Delta^n) \to C$ to the map $\varphi' : \on{spine}(\Delta^{n'}) \to C$ defined as follows: $\on{spine}(\Delta^{n'})$ is obtained by contracting each edge of $\on{spine}(\Delta^n)$ whose image under $\varphi$ is a degenerate edge of $C$, and $\varphi'$ is the unique factoring of $\varphi$ through the contraction $\on{spine}(\Delta^n) \to \on{spine}(\Delta^{n'})$. (The hypothesis in the lemma is needed in showing that this `right adjoint' is actually a functor.) Quillen's Theorem A implies that these functors are homotopy equivalences.
		
		If $a = b$, the 1-category $\ms{Spine}_C(a, a)$ splits as a disjoint union of 1-categories 
		\[
		\ms{Spine}_C(a, a) = \ms{Spine}_C^1(a, a) \sqcup \ms{Spine}_C^2(a, a)
		\]
		where the first factor consists of maps $\varphi : \on{spine}(\Delta^n) \to C$ sending each edge to the degenerate edge $\id_a$, and the second factor is its complement. The 1-category $\ms{Spine}_C^1(a, a)$ is contractible because it has an initial object $\on{spine}(\Delta^0) \to C$ sending the point to $a \in C$. (It also has a terminal object $\on{spine}(\Delta^1) \to C$ sending the edge to $\id_a$.) There is an adjunction 
		\[
		\ms{Spine}_C^{\mr{nd}}(a, a) \rightleftarrows \ms{Spine}_C^2(a, a) 
		\]
		defined as before, so the embedding $\ms{Spine}_C^{\mr{nd}}(a, a) \hra \ms{Spine}_C^2(a, a)$ is a homotopy equivalence. These two observations yield the desired statement. 	
	\end{myproof}
	
	\subsection{The reduced lift presentation} \label{ss-present} 
	
	\subsubsectiona \label{dob-def} 
	Interpret monoids in $\ms{Spaces}$ as $\infty$-categories with one object. Interpret groups in $\ms{Spaces}$ as $\infty$-groupoids with one object. Then we obtain functors
	\begin{cd}
		\left\{\parbox{1in}{\begin{center} Simplicial sets with one vertex \end{center}} \right\} \ar[r, "\text{Joyal}"] & \left\{\parbox{1in}{\begin{center} Monoids in $\ms{Spaces}$ \end{center}} \right\} \ar[r, "\text{Kan}"] & \left\{\parbox{1in}{\begin{center} Groups in $\ms{Spaces}$ \end{center}} \right\}
	\end{cd}
	The first functor is Joyal fibrant replacement, and the second functor is Kan fibrant replacement. These could be called `homotopy monoidification' and `homotopy groupification,' respectively. We will define a partial 1-category with one vertex (i.e.\ a \emph{partial monoid}) $S_I$ whose homotopy monoidification is $\BB_I^+$ (when $I$ is arbitrary) and whose homotopy groupification is $\BB_I$ (when $I$ is finite or affine type). 
	
	\begin{defn} \label{dob-present} 
		The \emph{finite type reduced lift presentation} $S_I$ is the partial 1-category with one vertex defined as follows: 
		\begin{itemize}
			\item There is an arrow $\alpha_w$ for each $w \in W_{I, \mr{fin}}$. 
			\item For $w_1, w_2 \in W_{I, \mr{fin}}$, the composite $\alpha_{w_1} \circ \alpha_{w_2}$ equals $\alpha_{w_1w_2}$ if $(w_1, w_2)$ is a reduced sequence, and it is undefined otherwise. 
		\end{itemize}
		The simplicial set $S_I$ has the following explicit characterization:  
		\begin{itemize}
			\item A simplex $\sigma : \Delta^{n_\sigma} \to S_I$ is a reduced sequence $(w_{\sigma, 1}, \ldots, w_{\sigma, n_\sigma})$ in $W_{I, \mr{fin}}$. 
			\item There is a diagram 
			\begin{cd}
				\Delta^{n_{\sigma_2}} \ar[rd, swap, "\sigma_2"] \ar[rr, "F"] & & \Delta^{n_{\sigma_1}} \ar[ld, "\sigma_1"] \\
				& S_I
			\end{cd}
			if and only if the following holds for all $0 < j \le n_{\sigma_2}$: 
			\[
			w_{\sigma_2, j} = w_{\sigma_1, F(j-1)+1}\  w_{\sigma_1, F(j-1)+2}\  \cdots \ w_{\sigma_1, F(j)}. 
			\]
		\end{itemize}
		Similarly, the \emph{reduced lift presentation} $S_I'$ is the partial 1-category obtained by reproducing the above definition with $W_{I, \mr{fin}}$ replaced by $W_I$. The maps 
		\[
			W_{I, \mr{fin}} \hra W_I \xra{r} \BB_I^+
		\]
		induce maps of partial 1-categories 
		\begin{equation} \tag{$\diamond$} \label{presenteq}
			S_I \to S_I' \to (\text{delooping of } \BB_I^+). 
		\end{equation}
	\end{defn}
	
	\begin{thm}{\cite[Thm.\ 5.1]{dob}} \label{dob-thm2}
		The maps in~\textnormal{(\ref{presenteq})} are categorical equivalences. 
	\end{thm}
	\begin{myproof}
		First, we show that $S_I \to \BB_I^+$ is a categorical equivalence. In view of Corollary~\ref{partial3} and Lemma~\ref{partial4}, it suffices to show that the map 
		\[
			p : \ms{Spine}_{S_I}^{\mr{nd}}(*,*) \to \BB_I^+
		\]
		is a homotopy equivalence. The classical `reduced lift' presentation of $\BB_I^+$ implies that this map induces a bijection on connected components. Therefore, we fix $b \in \BB_I^+$ and show that $p^{-1}(b)$ is contractible. By comparing Remark~\ref{results-rmk} with Definition~\ref{dob-present}, we obtain an equivalence of 1-categories 
		\[
			\ms{Word}_{\mr{fr}}(b) \simeq p^{-1}(b). 
		\]
		Corollary~\ref{results3} says that $\ms{Word}_{\mr{fr}}(b)$ is contractible, so $p^{-1}(b)$ is also contractible, as desired. 
		
		To show $S_I' \to \BB_I^+$ is a categorical equivalence, replace $\ms{Word}_{\mr{fr}}(b)$ by $\ms{Word}_{\mr{r}}(b)$ in the previous paragraph. 
	\end{myproof}
	
	\begin{cor} \label{cors1}
		We have 
		\[
		\BB_I^+ \simeq \colim_{\substack{J \subseteq I \\ J \textnormal{ finite type}}} \BB_J^+
		\]
		where the colimit is evaluated in the $\infty$-category of monoids in $\mathsf{Spaces}$.
	\end{cor}
	\begin{myproof}
		As in~\ref{dob-def}, we think of monoids in $\ms{Spaces}$ as $\infty$-categories with one object. By Theorem~\ref{dob-thm2}, we can model the delooping of $\BB_J^+$ by the simplicial set $S_J$. By~\cite[Thm.\ 4.2.4.1]{htt}, homotopy colimits (in model categories) compute $\infty$-categorical colimits, so it suffices to compute the homotopy colimit of the diagram $J \mapsto S_J$. Because monomorphisms are cofibrations in the Joyal model structure, the diagram $J \mapsto S_J$ is cofibrant in the projective model structure for the diagram 1-category (see~\cite[Prop.\ 5.12]{ozornova}), so the homotopy colimit $\colim_J S_J$ is equivalent to the colimit in the 1-category of simplicial sets. The latter is $S_I$, which models the delooping of $\BB_I^+$, as desired. 
	\end{myproof}
	
	\subsection{Application: intertwining functors} \label{ss-apply} 
	
	\begin{thm} \label{cors4} 
		Assume that $I$ is finite or affine type. The map from the delooping of $\BB_I^+$ to the delooping of $\BB_I$ is a homotopy equivalence. 
	\end{thm}
	\begin{myproof}
		Theorem~5.2 of \cite{dob} says that the homotopy groupification of $\BB_I^+$ is homotopy equivalent to the loop space of the $W_I$-quotient of the complexified root hyperplane arrangement corresponding to $I$. If $I$ is affine, the main theorem of~\cite{ps} says that the latter space is homotopy equivalent to $\BB_I$, and the claim follows. The analogous result when $I$ is finite type was established in~\cite{deligne}.
	\end{myproof}
	
	\begin{rmk*}
		The possibility of making this deduction was already noted in~\cite[Thm.\ 6.3]{dob}. At that time, \cite{ps} had not yet appeared, which is why~\cite[Cor.\ 6.5]{dob} restricts to finite type Coxeter systems.
	\end{rmk*}
	
	\begin{cor} \label{braid-group}
		Assume that $I$ is finite or affine type. Let $\mc{H}_I$ be a $\D(B\backslash G / B)$-version of the Hecke $\infty$-category with cells indexed by $W_I$.\footnote{See~\cite[6.1]{tatr} for a precise definition.} There is a monoidal functor $F : \BB_I \to \mc{H}_I$ which is determined essentially uniquely by the following requirements:
		\begin{enumerate}[label=(\roman*)]
			\item For any $w \in W_I$, we have $F(r(w)) = j_{w, *} \ul{\BC}$.\footnote{The right hand side is the $*$-extension of the constant rank-1 local system on the $w$ Schubert cell.} 
			\item For any reduced sequence $(w_1, w_2)$ in $W_I$, the monoidality isomorphism 
			\[
			F(r(w_1)) \otimes F(r(w_2)) \simeq F(r(w_1w_2))
			\]
			is induced by the identity automorphism of the constant rank-1 local system on the $w_1w_2$ Schubert cell.  
		\end{enumerate}
	\end{cor}
	\begin{myproof} 
		Giving a monoidal functor $F : \BB_I^+ \to \mc{H}_I$ is equivalent to giving a map of monoids in $\mathsf{Spaces}$ from $\BB_I^+$ to the maximal sub-$\infty$-groupoid of $\mc{H}_I$. By Theorem~\ref{dob-thm2}, this is equivalent to providing the following data: 
		\begin{enumerate}[label=(\roman*)]
			\item For any $w \in W_I$, an object $\mc{F}_w \in \mc{H}_I$.
			\item For any reduced sequence $(w_1, w_2)$ in $W_I$, an isomorphism $\mc{F}_{w_1} \otimes \mc{F}_{w_2} \simeq \mc{F}_{w_1w_2}$.
			\item For any reduced sequence in $W_I$ of length $\ge 3$, a higher compatibility between the isomorphisms chosen in (ii). 
		\end{enumerate}
		The first two pieces of data are specified by (i) and (ii) in the theorem statement. The higher compatibilities in (iii) live in the higher homotopy groups of the mapping spaces 
		\[
			\tau^{\le 0} \Hom_{\mc{H}_I}(j_{w, *}\ul{\BC}, j_{w, *}\ul{\BC})
		\]
		for various $w \in W_I$. (The Hom in a stable $\infty$-category is a spectrum, and the mapping space is $\Omega^\infty$ of that spectrum. Under the cohomological indexing convention, this corresponds to the $\tau^{\le 0}$ truncation.) These mapping spaces are discrete, so (iii) amounts simply to checking the 3-term associativity condition. This condition clearly holds for the choices of (i) and (ii) in the theorem.
		
		Next, we enhance this to a map from $\BB_I$. Each $\mc{F}_w:= j_{w, *}\ul{\BC}$ has a monoidal inverse given by $j_{w, !}\ul{\BC}[2\ell(w)]$. Hence, the map from $\BB_I^+$ extends to a map from the homotopy groupification of $\BB_I^+$. By Theorem~\ref{cors4}, the latter coincides with $\BB_I$.
	\end{myproof}

\end{document}